\documentclass[12pt]{article}
\usepackage{latexsym,amsfonts,amsthm,amsmath,amscd,amssymb}
\usepackage[dvips]{graphicx}

\newtheorem{lemma}{Lemma}[section]
\newtheorem{thm}[lemma]{Theorem}
\newtheorem{rem}[lemma]{Remark}
\newtheorem{prop}[lemma]{Proposition}
\newtheorem{cor}[lemma]{Corollary}

\newtheorem{defn}[lemma]{Definition}

\newcommand\matS{{\mathbb{S}}}
\newcommand\matP{{\mathbb{P}}}
\newcommand\matD{{\mathbb{D}}}

\newcommand\matR{{\mathbb{R}}}
\newcommand\matN{{\mathbb{N}}}

\newcommand\St{{St}}

\renewcommand{\hbar}{{\overline{h}}}

\newfont{\Got}{eufm10 scaled 1200}

\newcommand{\mycap} [1] {\caption{\footnotesize{#1}}}

\newcommand{\orb}{\mathop{\rm orb}\nolimits}

\newcommand{\Diff}{\mathop{\rm Diff}\nolimits}

\begin{document}

\title{Lower bounds on complexity of geometric $3$-orbifolds}

\author{Ekaterina~{\textsc Pervova}}

\maketitle

\begin{abstract}
We establish a lower bound on the complexity orientable locally
orientable geometric 3-orbifolds in terms of Delzant's T-invariants
of their orbifold-fundamental groups, generalizing previously known
bounds for complexity of 3-manifolds.

\end{abstract}

\section*{Introduction}

The aim of this note is to provide lower bounds on the so-called
\emph{complexity} of $3$-orbifolds \cite{CP01} in terms of
\emph{relative T-invariant} \cite{Del1,Del2} of their orbifold
fundamental groups. More precisely, we consider orientable locally
orientable geometric $3$-orbifolds and provide two types of bounds:
one in terms of the T-invariant relative to a finite set of
subgroups (which however is not uniquely defined) and another in
terms of the T-invariant relative to an elementary family of
(finite) subgroups, which is uniquely defined by the orbifold under
consideration. The bounds thus obtained generalize the bound on the
complexity \cite{Ma} of $3$-manifolds in terms of the T-invariants
of their fundamental groups, that was already established in
\cite{Del1}.

\bigskip

\noindent \textsc{Acknowledgements}: The author would like to thank
Thomas Delzant for useful discussions regarding the T-invariant of
groups.

\section{Main definitions}

In this section we briefly recall the notions of orbifold, its
complexity, and the $T$-invariant of a group.

\paragraph{Local structure of orbifolds} For the general theory of
orbifolds we refer the reader to \cite{thurston}, to the recent
\cite{BMP}, and to the extensive bibliography of the latter (the
essential definitions,as well as an overview of results concerning
geometrization of orbifolds, can also be found in the survey
\cite{Bon}). We just recall that an orbifold of dimension $n$ is a
topological space with a singular smooth structure, locally modelled
on a quotient of $\matR^n$ under the action of a finite group of
diffeomorphisms. We will confine ourselves to the case of compact
orientable locally orientable $3$-orbifolds. Such an orbifold $X$ is
given by a compact \emph{support} $3$-manifold $|X|$ together with a
\emph{singular set} $S(X)$. Given the assumption of local
orientability, $S(X)$ is a finite collection of circles and
univalent graphs tamely embedded in $|X|$, where the univalent
vertices are given by the intersection of $S(X)$ with $\partial|X|$
and each component of $S(X)$ minus the vertices is endowed with an
\emph{order} in $\{p\in\matN: p\geqslant 2\}$, with the restriction
that the three germs of edges incident to each vertex should have
orders $(2,2,p)$, for arbitrary $p$, or $(2,3,p)$, for
$p\in\{3,4,5\}$.

\paragraph{Uniformizable $n$-orbifolds} An \emph{orbifold-covering}
is a map between orbifolds locally modelled on a map of the form
$\matR^n/\Delta\rightarrow\matR^n/\Gamma$, naturally defined
whenever $\Delta<\Gamma<\Diff_+(\matR^n)$. If $M$ is a manifold and
$\Gamma$ is a group acting properly discontinuously on $M$ then
$M/\Gamma$ possesses a naturally defined structure of an orbifold
such that the natural projection $M\rightarrow M/\Gamma$ is an
orbifold-covering map. Any orbifold obtained as $M/\Gamma$ is said
to be \emph{uniformizable}. We note that in the case $n=2$ the set
of non-uniformizable $2$-orbifolds (also called \emph{bad}
$2$-orbifolds) admits the following easy description:

\begin{lemma}
The only non-uniformizable closed $2$-orbifolds are $(S^2;p)$, the
$2$-sphere with one cone point of order $p$, and $(S^2;p,q)$, the
$2$-sphere with cone points of orders $p\neq q$.
\end{lemma}

\paragraph{Spherical and discal orbifolds} We now introduce notation and terminology
for several particular types of $2$- and $3$-orbifolds (the notation
we employ follows \cite{CP01}). We define $\matD_o^3$ to be $D^3$,
the \emph{ordinary discal} $3$-orbifold, $\matD_c^3(p)$ to be $D^3$
with singular set a trivially embedded arc with arbitrary order $p$,
and $\matD_v^3(p,q,r)$ to be $D^3$ with singular set a trivially
embedded triod with edges of admissible orders $p$, $q$, $r$. We
call $\matD_c^3(p)$ and $\matD_v^3(p,q,r)$ respectively \emph{cyclic
discal} and \emph{vertex discal} $3$-orbifolds; we will occasionally
suppress the indication of orders and write $\matD_c^3$ and
$\matD_v^3$ to denote a generic cyclic or vertex discal
$3$-orbifold.

We also define the \emph{ordinary}, \emph{cyclic}, and \emph{vertex
spherical} $2$-orbifolds, denoted respectively by $\matS_o^2$,
$\matS_c^2(p)$, and $\matS_v^2(p,q,r)$, as the $2$-orbifolds
bounding the corresponding discal $3$-orbifolds $\matD_o^3$,
$\matD_c^3(p)$, and $\matD_v^3(p,q,r)$. Finally, we define the
\emph{ordinary}, \emph{cyclic}, and \emph{vertex spherical}
$3$-orbifolds, denoted respectively by $\matS_o^3$, $\matS_c^3(p)$,
and $\matS_v^3(p,q,r)$, as the $3$-orbifolds obtained by mirroring
the corresponding discal $3$-orbifolds $\matD_o^3$, $\matD_c^3(p)$,
and $\matD_v^3(p,q,r)$ in their boundary. We will also use the
notation $\matS_*^2$ and $\matS_*^3$ to indicate a spherical $2$- or
$3$-orbifold with generic orders.

\paragraph{$2$-suborbifolds and irreducible $3$-orbifolds} We say
that a $2$-orbifold $\Sigma$ is a \emph{suborbifold} of a
$3$-orbifold $X$ if $|\Sigma|$ is embedded in $|X|$ so that
$|\Sigma|$ meets $S(X)$ transversely (in particular, it does not
meet the vertices), and $S(\Sigma)$ is given precisely by
$|\Sigma|\cap S(X)$, with matching orders.

A spherical $2$-suborbifold $\Sigma$ of a $3$-orbifold $X$ is called
\emph{essential} if it does not bound in $X$ a discal $3$-orbifold.
A $3$-orbifold $X$ is called \emph{irreducible} if it does not
contain any non-uniformizable $2$-suborbifold and every spherical
$2$-suborbifold of $X$ is inessential (in particular, it is
separating).

\paragraph{Geometric $3$-orbifolds} An orbifold is said to be
\emph{geometric} if it admits a complete geometric structure in the
sense of orbifolds. Namely, let $\tilde{X}$ be a (simply) connected
manifold, and let $G$ be a group acting effectively and transitively
on $X$ such that all point stabilizers are compact. An orbifold
admits a geometric structure modelled on $(\tilde{X},G)$ if its
maximal orbifold atlas contains an atlas such that all the folding
groups consist of appropriate restrictions of the elements of $G$
and all the charts are again restrictions of elements of $G$.

In particular, it can be shown that any geometric $3$-orbifold $X$
is isomorphic to an orbifold of the form $\tilde{X}/\Gamma$, where
$\tilde{X}$ is a simply connected manifold and $\Gamma$ acts
properly discontinuously on $\tilde{X}$, see for instance
\cite{thurston}. The group $\Gamma$ (viewed up to isomorphism) is
called the \emph{orbifold-fundamental group} of $X$ (a more direct
definition can be found, for instance, in \cite{Scott}).

We note, although we won't need this, that the analogue of
geometrization conjecture for orbifolds was announced by Thurston in
1982 and that the complete proof of the conjecture can be found in
\cite{Boileau}.

\begin{figure}
    \begin{center}
    \includegraphics[scale=0.5]{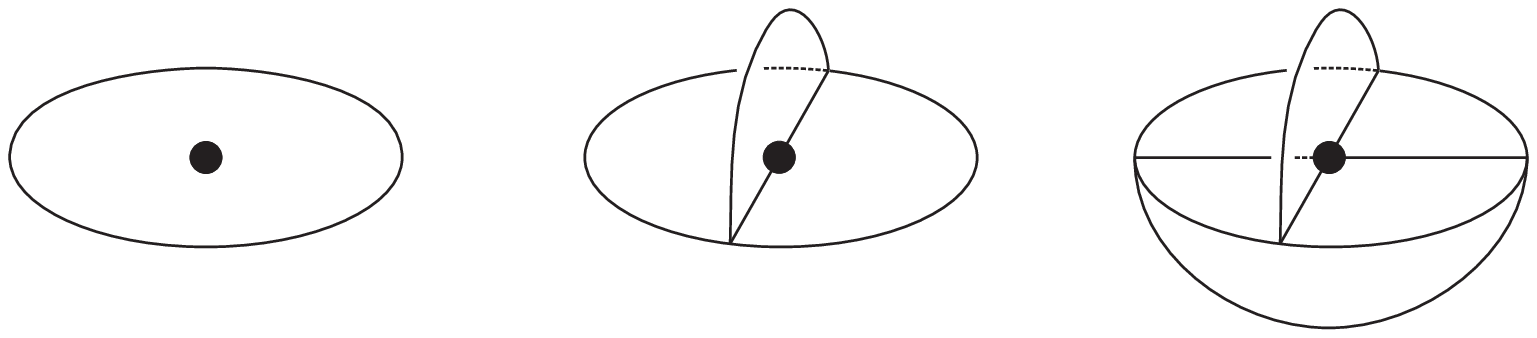}
    \mycap{Neighbourhoods of points in a simple polyhedron}
    \label{neighborhds_pts_simple:fig}
    \end{center}
    \end{figure}

\paragraph{Simple and special polyhedra} From now on we will employ
the piecewise linear viewpoint, which is equivalent to the smooth
one in dimensions $2$ and $3$. A \emph{simple polyhedron} is a
compact polyhedron $P$ such that the link of each point of $P$ can
be embedded in the space given by a circle with three radii. In
particular, $P$ has dimension at most $2$. A point of a simple
polyhedron is called a \emph{vertex} if its link is given precisely
by a circle with three radii; a regular neighbourhood of a vertex is
shown in Fig.~\ref{neighborhds_pts_simple:fig} on the far right. It
follows from the figure that the vertices are isolated, hence finite
in number (the case when there are no vertices is also allowed). The
complexity $c(P)$ of $P$ is the number of vertices that $P$
contains.

Two more restrictive types of polyhedra will be used below. A simple
polyhedron $P$ is called \emph{almost-special} if the link of each
point of $P$ is given by a circle with zero or two or three radii.
The three possible types of neighbourhoods of points of $P$ are
shown in Fig.~\ref{neighborhds_pts_simple:fig}. The points of the
second and the third type are called \emph{singular}, and the set of
all singular points of $P$ is denoted by $S(P)$. We say that $P$ is
special if it is almost-special, $S(P)$ contains no circle
component, and $P\setminus S(P)$ consists of open $2$-discs.

\paragraph{Spines of $3$-orbifolds and complexity} The notions we
define in this paragraph were introduced in \cite{CP01} and
generalize the corresponding notions for $3$-manifolds. Let $X$ be a
closed orientable locally orientable $3$-orbifold, and let $P$ be a
simple polyhedron contained in $|X|$. We say that $P$ is a
\emph{spine} of $X$ if the following conditions are satisfied:
\begin{itemize}
  \item The intersections between $P$ and $S(X)$ occur only in
  surface points of $P$ and non-vertex points of $S(X)$, and they
  are transverse;
  \item If $U(P)$ is an open regular neighbourhood of $P$ in $|X|$
  then each component of $X\setminus U(P)$ is isomorphic to one of
  the discal $3$-orbifolds $\matD_*^3$.
\end{itemize}
We will normally identify $X\setminus U(P)$ with $X\setminus P$,
making a distinction only when it matters; it is easy to check that
this choice does not cause any ambiguity.

If $P$ is a spine of $X$ as above, we define the following function
(that depends not only on $P$ as an abstract polyhedron but also on
its position in $|X|$ relative to $S(X)$):
$$c(P,S(X))=c(P)+\sum\{p-1\,:\, x\in P\cap S(X),\mbox{ the order of }S(X)\mbox{ at }x\mbox{ is }p\}.$$
Then the \emph{complexity} $c(X)$ of the $3$-orbifold $X$ is defined
as the minimum of $c(P,S(X))$ taken over all spines $P$ of $X$. We
note that $c$ is always well-defined because every orbifold has
simple spines \cite[Section 2]{CP01}.

We say that a spine $P$ of $X$ is \emph{minimal} if $c(P,S(X))=c(X)$
and every proper subpolyhedron of $P$ which is also a spine of $X$
is actually homeomorphic to $P$. We cite the following property of
such spines (actually established in \cite{CP01} under some slightly
relaxed assumptions), that will be important for us:

\begin{lemma}\label{cpt:once:S(X):lem}\emph{\cite[Lemma 2.1]{CP01}}
Let $P$ be a minimal spine of an orbifold $X$, and let $\alpha$ be a
connected component of $S(X)$ minus the vertices. Then $\alpha\cap
P$ consists of precisely one point.
\end{lemma}

\paragraph{Orbifolds with minimal special spines} In addition to the
orbifold $\matS_*^3$ described above, we also need to consider
certain orbifolds $(\matP^3,F_p)$ and $(L_{3,1},F_p)$. In both cases
$F_p$ is a circle of order $p$, given by a non-singular fibre of the
natural Seifert fibration; since we will have to include the case of
manifolds $S^3$, $\matP^3$, $L_{3,1}$, we stipulate that if $K$ is a
knot then $K_p$ denotes $K$ equipped with the cone order $p$ if
$p\geqslant 2$, and it denotes the empty set if $p=1$. With this
stipulation we can write $\matS_o^3=\matS_c^3(1)$.

In the above notation, it was established in \cite[Proposition
2.4]{CP01} that the complexity of any of $\matS_c^3(p)$,
$(\matP^3,F_p)$, $(L_{3,1},F_p)$ is equal to $p-1$; moreover, any
minimal spine of any of these orbifolds is not special.

On the other hand, we have the following result, also established in
\cite{CP01}.

\begin{thm}\label{excep:orbs:thm} \emph{\cite[Theorem 2.6]{CP01}}
If $X$ is an irreducible $3$-orbifold and not $\matS_c^3(p)$,
$(\matP^3,F_p)$, $(L_{3,1},F_p)$, or $\matS_v^3(p,q,r)$, then any
minimal spine of $X$ is special.
\end{thm}

\paragraph{Duality} We define a \emph{triangulation} of a
$3$-orbifold $X$ to be a triangulation of the manifold $|X|$ such
that $S(X)$ is a subset of the $1$-skeleton. Note that by a
triangulation of a $3$-manifold $M$ we mean a realization of $M$ as
a gluing of a finite number of tetrahedra along a complete system of
simplicial pairings of the lateral faces and that, in particular, we
allow multiple and self-adjacencies of the tetrahedra. Recall now
that there is a well-known duality between the set of such
triangulations of $M$ and the set of special polyhedra $P$ embedded
in $M$ in such a way that $M\setminus P$ is a union of open
$3$-discs; indeed, any $P$ as above induces a natural
cellularization of $M$ the dual of which is a triangulation, and
vice versa. This duality extends, with reservations, also to the
case of $3$-orbifolds.

\begin{prop}\label{dual:tri:orb:prop}\emph{\cite[Proposition 2.8]{CP01}}
Let $X$ be a $3$-orbifold as in Theorem \ref{excep:orbs:thm}. Then
dual to any minimal spine of $X$ there is a triangulation (in the
sense of orbifolds).
\end{prop}

\paragraph{T-invariant relative a set of subgroups } Finally, we need
the following definition introduced in \cite{Del1}.

\begin{defn}\label{relative_T-invariant}
Let $G$ be a group, and let $(C_1,\ldots,C_n)$ be a family of its
subgroups. Then the {\em T-invariant} of the pair
$(G;C_1,\ldots,C_n)$ is defined by the following condition:
$T(G;C_1,\ldots,C_n)\leqslant k$ if and only if there exists a
simply connected simplicial polyhedron $\Pi$ of dimension $2$ such
that $G$ acts on $\Pi$, this action is simplicial and possesses the
following properties:
\begin{itemize}
\item the number of faces of $\Pi$ modulo the action of $G$ is no
more than $k$;

\item the stabilizers of vertices of $\Pi$ are conjugate to $C_1$,
$\ldots$, $C_n$ and each of these groups fixes some vertex of
$\Pi$.
\end{itemize}
\end{defn}

\section{Lower bounds}

In this section we establish two types of lower bounds on the
complexity of $3$-orbifolds that are mentioned in the introduction.

\paragraph{Bounds with respect to a finite family}
Consider a geometric orbifold $X$. Denote by $\tilde{X}$ the
universal orbifold covering of $X$ and by $\pi$ the covering of $X$
by $\tilde{X}$ induced by the action of $\pi_1^{\orb}(X)$ on
$\tilde{X}$. To proceed, we recall the following easy facts:

\begin{itemize}
\item Let $x\in\tilde{X}$ be a point such that $\pi(x)\in S(X)$,
the singular set of $X$. Then the stabilizer $\St(x)$ of $x$ in
$\pi_1^{\orb}(X)$ is a finite subgroup.

\item For every two points $x,y\in\tilde{X}$ such that $\pi(x)$,
$\pi(y)$ belong to the same connected component of $S(X)$ minus the
vertices, the orders of subgroups $\St(x)$ and $\St(y)$ are equal to
each other and to $p$, the cyclic order associated to that
component.
\end{itemize}

Denote by $c_X$ the number of circular components in $S(X)$ minus
vertices and by $v_X$ the number of vertices in $S(X)$.

\begin{thm}\label{lower:bound:thm}
Let $X$ be a closed orientable locally orientable geometric
$3$-orbifold with non-empty singular set such that $X$ has a minimal
special spine. Then there is a finite family of subgroups
$C_1,\ldots,C_{c_X+v_X}\leqslant \pi_1^{\orb}(X)$ such that
$$c(X)\geqslant \frac{1}{2}T(\pi_1^{\orb}(X);C_1,\ldots,C_{c_X+v_X})+\sum_{\alpha\subset S(X)}(p-1), $$
where the latter sum is taken over all connected components
$\alpha$ of $S(X)$ minus the vertices, and $p$ is the order
associated to $\alpha$.
\end{thm}

\begin{rem}
\emph{As follows from Theorem \ref{excep:orbs:thm}, most geometric
$3$-orbifolds admit minimal special spines.}
\end{rem}

\begin{proof}
Let $P$ be a minimal special spine of $X$. It follows from
Lemma~\ref{cpt:once:S(X):lem} that
$$c(X)=\#V(P)+\sum_{\alpha\subset S(X)}(p-1),$$
where $\alpha$, $p$ are as in the statement of the theorem. So we
actually need to estimate $\#V(P)$.

Let $\tau$ be any triangulation of $X$ defined in
Proposition~\ref{dual:tri:orb:prop} (\emph{i.e.}, dual to $P$ and
such that $S(X)$ is contained in the 1-skeleton of $\tau$). This
triangulation lifts to a triangulation $\tilde{\tau}$ of
$\tilde{X}$. Denote by $Q$ the 2-skeleton of $\tau$ and by
$\tilde{Q}$ the 2-skeleton of $\tilde{\tau}$. Then obviously
$\pi_1^{\orb}(X)$ acts on $\tilde{Q}$, and
$\tilde{Q}/\pi_1^{\orb}(X)=Q$.

Since $S(X)$ is contained in the 1-skeleton of $\tau$ and
therefore of $Q$, its pre-image $\pi^{-1}(S(X))$ is contained in
the 1-skeleton of $\tilde{Q}$. Hence the action of
$\pi_1^{\orb}(X)$ is simplicial, and all points with nontrivial
stabilizers are contained in the 1-skeleton of $\tilde{Q}$.

Notice that among the stabilizers of vertices of $\tilde{Q}$ there
is only a finite number of non-conjugate subgroups. Indeed, if
$\pi(u)=\pi(v)$ then $\St(u)$ is conjugate to $\St(v)$, and $Q$
contains only a finite number of vertices. We claim that, if $S(X)$
is non-empty then the number $\#V(Q)$ of vertices of $Q$ is
$c_X+v_X$.

Since triangulation $\tau$ is dual to $P$, this number is exactly
the number of connected components of $X\setminus P$. Obviously,
this number is at least $c_X+v_X$ because each vertex of $S(X)$ must
lie in some ball of $X\setminus P$, and that ball cannot contain any
other vertices of $S(X)$ or have non-empty intersection with any
circular components. Also, every circular component must intersect
some ball, and it will be the only connected component of $S(X)$
intersecting that ball.

On the other hand, suppose that $S(X)$ is non-empty and the
complement $X\setminus P$ contains an ordinary discal 3-orbifold
$B$. Then there is a 2-cell $c$ of $P$ separating $B$ from another
discal 3-orbifold $B'$. Obviously, $c$ does not contain any points
of $S(X)$. So we can puncture $c$ fusing $B$ and $B'$ into a discal
3-orbifold and obtain a new simple spine $P'$ of $X$. Notice that
$c$ must contain some vertices, so $\#V(P')<\#V(P)$ and hence
$c(P',S(X))<c(P,S(X))$. This contradicts the minimality of $P$.
Hence every component of $X\setminus P$ has non-empty singular set.

Finally, by Lemma~\ref{cpt:once:S(X):lem} each circular component
$\alpha$ of $S(X)$ is contained in the closure of precisely one
component of $X\setminus P$, which approaches the cell intersected
by $\alpha$ from both sides. By the same lemma each non-closed edge
of $S(X)$ intersects exactly two (with multiplicity) components of
$X\setminus P$, namely the ones containing its ends. This and the
absence of ordinary discal 3-orbifolds implies that any component of
$X\setminus P$ contains either a vertex of $S(X)$ or a portion of a
circular component. Hence the number of components of $X\setminus P$
is indeed $c_X+v_X$.

In particular, the set of vertices of $\tilde{Q}$ is decomposed
into $c_X+v_X$ classes of vertices such that two vertices are in
the same class if and only if they project to the same vertex of
$Q$. The stabilizers of vertices in the same class are conjugate,
and we choose exactly one of them. Denote the chosen subgroups by
$C_1$, $\ldots$, $C_{c_X+v_X}$. Notice that $\tilde{Q}$ and
subgroups $C_1$, $\ldots$, $C_{c_X+v_X}$ satisfy the second
condition in the definition of the T-invariant.

Since $\tilde{X}\setminus\tilde{Q}$ consists of open balls,
$\tilde{Q}$ is simply-connected. Therefore we have that
$$\#\mbox{2-faces of }Q\geqslant T(\pi_1^{\orb}(X);C_1,\ldots,C_{c_X+v_X}).$$

Since $\tau$ is dual to $P$, we have that the number of 2-faces of
$P$, and therefore of $Q$, is equal to the number of edges of $P$.
Since $P$ is special, its singular graph has valence 4. Therefore
the number of edges of $P$ is twice the number of vertices of $P$.
Summarizing, we have:
$$\#V(P)\geqslant \frac{1}{2}T(\pi_1^{\orb}(X);C_1,\ldots,C_{c_X+v_X}).$$
This proves the theorem.
\end{proof}

\paragraph{Bounds with respect to an elementary family of subgroups}
The subgroups $C_1$, $\ldots$, $C_{c_X+v_X}$ that appear in the
statement of Theorem~\ref{lower:bound:thm} are not well-defined by
$X$ only; for instance, those among them that are stabilizers of
points projecting to vertices of $S(X)$ are well-defined up to
conjugacy only. However, the set of the orders of all these
subgroups are uniquely defined by the orbifold $X$. Thus, we can
obtain a lower bound on the orbifold complexity in terms depending
only on $X$, by employing the following definition introduced
in~\cite{Del3}.

\begin{defn}
Let $G$ be a fixed group. A family $\mathcal{C}$ of subgroups of
$G$ is called {\em elementary} if the following four conditions
hold:
\begin{enumerate}
\item[(a)] if $C\in\mathcal{C}$ then any subgroup of $C$ and any
subgroup conjugated to $C$ is in $\mathcal{C}$;

\item[(b)] any infinite subgroup in $\mathcal{C}$ is included into a
unique maximal subgroup in $\mathcal{C}$, and any increasing union
of finite subgroups in $\mathcal{C}$ is in $\mathcal{C}$;

\item[(c)] if $C\in\mathcal{C}$ acts on a tree then it fixes a point of
that tree or a point of its boundary, or it preserves a pair of
points of the boundary possibly permuting them;

\item[(d)] if $C\in\mathcal{C}$ is infinite, maximal in $\mathcal{C}$,
and $gCg^{-1}=C$ then $g\in C$.
\end{enumerate}
\end{defn}

In particular, the family of all finite subgroups of order not
greater than $p$, where $p$ is a fixed number, is an elementary
family. Denote it by $\mathcal{C}_p(G)$. We have the following
definition, from~\cite{Del3} as well.

\begin{defn}
The \emph{T-invariant} $T_p(G)$ relative to $\mathcal{C}_p(G)$ of a
group $G$ is defined as $\inf T(G;C_1,\ldots,C_n)$, where the
infimum is taken over all finite sets of subgroups from
$\mathcal{C}_p(G)$.
\end{defn}

\begin{rem}
\emph{The definition given in \cite{Del3} actually differs from the
one given above, in the following respects:
\begin{enumerate}
\item[1)] In the former, it is required that $G$ act on $\Pi$ from
Definition~\ref{relative_T-invariant} without edge inversions, and
that
\item[2)] the stabilizers of all edges and faces also be from
$\mathcal{C}_p$.
\end{enumerate}
However, we observe that in our situation (\emph{i.e.} that of
orbifold-fundamental groups of orientable locally orientable
$3$-orbifolds) these requirements are automatically satisfied.
Indeed, the stabilizer of each edge $\alpha$ is contained in the
intersection of stabilizers of its interior points. All those points
project to the same connected component of $S(X)\setminus V(S(X))$,
the set of vertices. Hence the orders of their stabilizers are all
equal to $p$, the local order associated to that connected
component. Hence the order of the stabilizer of $\alpha$ is no more
than $p$, so it belongs to the family $\mathcal{C}_{p(X)}$. Also,
since $\pi^{-1}(S(X))$ is contained in the 1-skeleton of
$\tilde{Q}$, all faces have trivial stabilizers. Finally, since $X$
is locally orientable, there are no inversions of edges.}
\end{rem}

Denote now by $p(X)$, where $X$ is an orbifold, the maximum of its
local orders. We now have the following result.

\begin{cor}
Let $X$ be an orbifold as in Theorem~\ref{lower:bound:thm}. Then
$$c(X)\geqslant \frac{1}{2}T_{p(X)}(\pi_1^{\orb}(X))+\sum_{\alpha\subset S(X)}(p-1).$$
\end{cor}

\vspace{1cm}

\noindent Centro De Giorgi (Scuola Normale Superiore)\\
Piazza dei Cavalieri, 3\\
56100 PISA -- Italy\\
\ \\
pervova@guest.dma.unipi.it\\

\end{document}